\documentclass[reqno,11pt]{amsart}
\usepackage{geometry}
\usepackage[all]{xy}
\usepackage{mathtools,amssymb,amsthm,mathrsfs,color,lineno,paralist,graphicx,float}
\usepackage[colorlinks,
linkcolor=red,
anchorcolor=green,
citecolor=blue,
]{hyperref}

\usepackage[T1]{fontenc}
\usepackage[utf8]{inputenc}
\usepackage[english]{babel}
\usepackage{graphicx}
\usepackage{subfigure}
\usepackage{tikz}
\usepackage{enumitem}

\usepackage{calc}
\definecolor{bleu1}{RGB}{0,57,128}
\def\bleu1{\color{bleu1}}

\usepackage{etoolbox}
\patchcmd{\section}{\normalfont}{\normalfont \bleu1}{}{}
\patchcmd{\subsection}{\normalfont}{\normalfont \bleu1}{}{}
\patchcmd{\subsubsection}{\normalfont}{\normalfont \bleu1}{}{}

\newcommand{\R}{{ \mathbb R}}

\newcommand{\T}{{ \mathbb T}}

\newcommand{\supp}{\operatorname{supp}}

\newtheorem{theorem}{Theorem}
\newtheorem{corollary}[theorem]{Corollary}

\newtheorem*{mainthm}{Main Theorem}

\theoremstyle{definition}
\newtheorem{remark}[theorem]{Remark}

\newtheorem*{definition*}{Definition}

\title{On the rigidity of the stable norm and Mather's $\beta$-function for geodesic flows}
\author{Anna Florio}
\address{CEREMADE, Université Paris-Dauphine, Université PSL, CNRS, 75016 Paris, France}
\email{florio@ceremade.dauphine.fr}

\author{Martin Leguil}
\address{École polytechnique, CMLS, Route de Saclay, 91128 Palaiseau Cedex, France}
\email{martin.leguil@polytechnique.edu}

\author{Alfonso Sorrentino}
\address{Dipartimento di Matematica, Università degli Studi di Roma ``Tor Vergata'', Via della ricerca scientifica 1, 00133 Rome, Italy}
\email{sorrentino@mat.uniroma2.it}

\date{\today}

\begin{document}

\maketitle
\begin{abstract} {We investigate rigidity phenomena associated to the stable norm and Mather's $\beta$-function for Riemannian geodesic flows on closed manifolds. Given two metrics $g_1$ and $g_2$, we compare these objects pointwise at individual homology classes. Our main result establishes that if Mather's $\beta$-function (or the stable norm) of $g_2$ at a non-zero homology class $h$ equals that of $g_1$ at $h$ multiplied by a suitable factor determined by the metrics, then the two metrics are homothetic on the Mather set of homology $h$ associated to $g_1$. In the case of conformally equivalent metrics, this yields a pointwise criterion for  homothety on the projected Mather set. Some consequences are discussed, including a pointwise rigidity result on the $2$-torus implying that if a metric has the same Mather's $\beta$-function at some non-zero homology class as a normalized  flat metric in the same conformal class, then the metric must be flat. This result can be considered a pointwise version of a similar global result by Bangert. Finally, an extension of these results to Ma\~n\'e's perturbations of general Tonelli Lagrangians is discussed.}
\end{abstract}

\section{Introduction}

The study of Riemannian geodesic flows reveals several fundamental objects that capture the asymptotic behavior of minimal geodesics and play a crucial role in understanding their structural and dynamical properties. Two key objects central to this work---the~\emph{stable norm} and the \emph{Mather’s $\beta$-function}---can be interpreted as manifestations of a homogenization effect at large scales.\\

 Let $(M, g)$ be a closed Riemannian manifold. We denote by $\|\cdot\|_x$ the associated norm in each fiber $T_xM$. Let us begin by recalling the definitions of stable norm and Mather's $\beta$-function and their relationship.\\

The \textit{stable norm} provides a geometric way to measure real homology classes by capturing the large-scale, or homogenized, geometry of $(M, g)$:
on a macroscopic level, the Riemannian metric ``behaves'' like  a norm on the homology group of the manifold. 
 It was introduced in the context of geometric measure theory by Federer~\cite{Federer74} and later became a staple of asymptotic Riemannian geometry through the work of Gromov~\cite{GromovLafontainePansu81} and others~\cite{BuragoBuragoIvanov01}.
More precisely, let $H_1(M; \mathbb{R})$ denote the first real homology group of $M$. The stable norm $\|\cdot\|_{g,\, \text{stable}}: H_1(M; \mathbb{R})\to \mathbb{R}$ is a function
\begin{equation*}
    \|h\|_{g,\,\text{stable}} := \inf \left\{ \sum |r_i| \, \mathcal{L}_g(\gamma_i) \right\}\,,
\end{equation*}
where the infimum is taken over all cycles $\sum r_i \gamma_i$ representing the homology class $h$, with $\gamma_i$ being $1$-simplices, and where $\mathcal{L}_g(\gamma_i)$ is the $g$-length of $\gamma_i$. Although the infimum may not be achieved in general,  when
the dimension of $M$ is two, it is indeed attained for every integer homology class (see \cite[Proposition 2.1]{CMP} or \cite[Proposition 5.6]{balacheff_massart_2008}).
This function indeed defines a norm on $H_1(M; \mathbb{R})$; we refer to \cite{GromovLafontainePansu81, BuragoBuragoIvanov01} for a comprehensive treatment.\\

From a dynamical perspective, one can consider the geodesic flow on $(M,g)$. It is the Euler-Lagrange flow of the Lagrangian $L_g\colon TM \to \R$ given by
\begin{equation*}
    L_g(x, v) = \frac{1}{2} \|v\|_x^2.
\end{equation*}
Let $\mathcal{M}(TM)$ denote the set of Borel probability measures on $TM$ that are {\it closed} and have finite first momentum.
We recall that a Borel probability measure $\mu$ on $TM$ is said to be closed if $\int_{TM} df_x(v)\, d\mu =0$ for every $f\in C^1(M)$, while it has 
finite first momentum if $\int_{TM}\|v\|_x \,d\mu< \infty$. We endow this space with the topology for which $\lim_{n \to +\infty }\mu_n =\mu$ if and only if
$\lim_{n\to \infty} \int_{TM} f\, d\mu_n =\int_{TM} d\, d\mu$ for any continuous function  $f\colon TM \rightarrow \R$ with linear growth, i.e.,
$$
\sup_{(x,v)\in TM} \frac{|f(x,v)|}{1+ \|v\|_x} < +\infty.
$$

\smallskip

To any $\mu \in \mathcal{M}(TM)$, one can associate its \textit{rotation vector} (or \textit{Schwartzman asymptotic cycle}) $\rho(\mu) \in H_1(M;\R)$. The resulting map $\rho\colon\mathcal{M}(TM) \to H_1(M;\R)$ is surjective and continuous (see \cite[Proposition 3.2.2]{SorrLectNotes}).
The {\it Mather's $\beta$-function} is then defined as:
\begin{equation*}
    \begin{aligned}
\beta_{g}\colon H_1(M;\R) &\longrightarrow \R\\
h &\longmapsto \beta_g(h):= \min_{\mu\in \rho^{-1}(h)} \int_{TM} L_g(x,v)\, d\mu.
    \end{aligned}
\end{equation*}

The probability measures achieving this minimum are invariant and are called \textit{Mather measures (or minimizing measures) of rotation vector $h$}. Their collection is denoted by $\mathfrak{M}_{g}^h$. The union of their supports defines an invariant set, called the \textit{Mather set of homology class $h$}:
\begin{equation}\label{def:mather_set}
    \widetilde{\mathcal{M}}_{g}^h = \bigcup_{\mu \in \mathfrak{M}_{g}^h} \supp(\mu) \subset TM.
\end{equation}
A cornerstone of the theory is Mather's graph theorem, which asserts that the projection $\pi\colon TM \to M$ restricts to a bi-Lipschitz homeomorphism from $\widetilde{\mathcal{M}}_{g}^h$ onto its image $\pi(\widetilde{\mathcal{M}}_{g}^h) \subset M$ (see \cite[Theorem 2]{Mather91} or \cite[Theorem 3.2.7]{SorrLectNotes}).
We call $\mathcal{M}_g^h := \pi(\widetilde{\mathcal{M}}_{g}^h)$ the {\it projected Mather set}.
\\

\bigskip

A fundamental result, which bridges the above geometric and dynamical viewpoints, is the following relation \cite[Proposition 1.4.2]{Massart96} (see also \cite{massart_1997}):
\begin{equation}\label{eq:beta_stable_relation}
    \beta_g(h) = \frac{1}{2} \|h\|_{g,\,\text{stable}}^2 \qquad \forall\, h \in H_1(M; \mathbb{R}).
\end{equation}
Thus, the stable norm and Mather's $\beta$-function are dual manifestations of the same underlying asymptotic geometry.\\

A natural  question arises: {\it to what extent do these objects allow one to recover information on the underlying Riemannian metric $g$?}\\

 This question can be approached from  different directions:

\begin{enumerate}[label=\textbf{(Q\arabic*)},ref=(Q\arabic*)]
    \item\label{question_un}  
    How do the local or global regularity properties of $\|\cdot\|_{g,\,  \text{stable}}$ or $\beta_g$ reflect the fine structure of the geodesic flow?

    	\smallskip
    \item\label{question_deux} 
      Can a metric  or a class of metrics be identified from knowledge of their stable norms or Mather's $\beta$-functions?
\end{enumerate}

\medskip

The differentiability properties of Mather's $\beta$-function are subtle and, in some cases, intimately tied to the underlying dynamics of the geodesic flow and the structure of Mather sets. Fine differentiability properties of $\beta_g$ have been thoroughly investigated in the case of closed  surfaces, whereas the situation in higher dimensions presents greater
complexity and defies simple characterization; we refer to \cite{balacheff_massart_2008, Massart96, massart_1997, massart_2003, massart_2009} for several interesting results in this direction. \\

A very interesting answer  to~\ref{question_un} was provided in \cite[Theorem 5.3]{Bangert} (see also \cite{Bangert1, Mather90, massart_sorrentino_2011} for related results). Recall that a homology class $h$ is \textit{$k$-irrational} if $k$ is the dimension of the smallest subspace of $H_1(M,\mathbb{R})$ generated by integer classes and containing $h$. 

\medskip 
\begin{theorem}[Bangert]
\noindent {\it  When $M=\T^2$, $\beta_g$ is differentiable at a $1$-irrational 
	homology class if and only if the Mather set $\widetilde{\mathcal{M}}_{g}^h$ consists of an invariant  torus foliated by periodic orbits}.
\end{theorem}
\medskip

On the other hand, question~\ref{question_deux} is more delicate, since these objects are built from~\textit{minimizing} geodesics and measures, hence they typically only capture information on the subset of the phase space traversed by such curves. 

A positive result in this direction can be found in \cite[Theorem 6.1]{Bangert}: 

\medskip
\begin{theorem}[Bangert]\label{thm_bangert}
	Suppose a Riemannian metric $g$ on the $2$-torus $\T^2$ has the same stable norm or Mather's $\beta$-function as a flat metric $g_0$ on $\T^2$. Then, $(\T^2, g)$ and $(\T^2, g_0)$ are isometric by an isometry homotopic to the identity.
\end{theorem}

\medskip
The proof relies on Hopf's theorem, which states that a $2$-torus without conjugate points is flat. In particular, the absence of conjugate points  follows from the global differentiability
of Mather's $\beta$-functions and the above mentioned consequence of the differentiability at $1$-irrational homology classes.

A similar rigidity result was proved by \cite{Osuna} on the $n$-dimensional torus, $n\geq 3$, but requires that not only the stable norm of the metric coincides with that of a flat one, but also its  $(n-1)$-dimensional counterpart (i.e., the one obtained by a similar construction on $H_{n-1}(M;\R))$.

\medskip

\begin{remark}
We stress that without the assumption on the absence of conjugate points, the stable norm and Mather's $\beta$-function are very weak invariants. For example,
on a Riemannian $2$-torus ``with a big bump'' (see \cite[ p.~46]{Bangert1}), there exists an open set on which one can arbitrarily increase the metric without changing the stable norm.\\
This construction can be generalized starting from any $(M,g)$ satisfying
\begin{equation}\label{subset}
\Sigma_g:=\cup_{h\in H_1(M;\R)} {\mathcal M}^h_g \subsetneq M.
\end{equation}
One considers any smooth function $\lambda\ge 0$  on  $M$, not identically zero, that has  support in the complement of the closure of ${\Sigma_g}$. Then, one checks that the metrics $g_1:= e^\lambda g$ and $g$ have the same  Mather's $\beta$-functions.
\end{remark}

\bigskip

Therefore, a more reasonable question is the following:
\begin{enumerate}[label=\textbf{(Q\arabic*)},ref=(Q\arabic*)]\setcounter{enumi}{2}
    \item\label{question_quattre} For a closed Riemannian manifold $(M,g)$, does the stable norm or Mather's $\beta$-function {\it locally}  determine the metric $g$ on the Mather sets? For example, 
can specific features of the stable norm or Mather's $\beta$-function reveal whether two metrics are related on their Mather sets? 
    \end{enumerate}
     \medskip

 \medskip   
 In the present paper, we  explore  question~\ref{question_quattre}. More specifically, we conduct a pointwise comparison at individual homology classes, inspired by analogous results for billiard  systems in  \cite{BarBiaSor}, and deduce, among other things, a pointwise version of the above-mentioned Bangert's Theorem \ref{thm_bangert} (see Corollary \ref{maincordim2}). Moreover, in Section \ref{manelagrangian}, an extension of these results to Ma\~n\'e's perturbations of general Tonelli Lagrangians is discussed (see Theorem \ref{thmmane}).

\medskip

\begin{remark} 
On the $2$-torus, the $\beta$-function is closely related to the so-called \emph{marked length spectrum}. 
 In~\cite{ColLefPat}, for Riemannian surfaces with Anosov geodesic flow, the authors show that the marked length spectrum determines the conformal class of the underlying complex structure in Teichmüller space (see \cite[Proposition 3.7]{ColLefPat}). By a result of Katok~\cite{Katok}, this in turn implies that the marked length spectrum actually determines the metric up to isometry.\\
By the uniformization theorem, every Riemannian metric on  $\T^2$ is conformally equivalent to a flat metric. We may then ask whether the stable norm or Mather's $\beta$-function determines the class
	of the underlying complex structure in the Teichmüller space of $\T^2$. That is:
	\begin{enumerate}[label=\textbf{(Q\arabic*)},ref=(Q\arabic*)]\setcounter{enumi}{3}
    \item\label{question_teich}	 Given two metrics $g_1$ and $g_2$ on $\T^2$ with the same Mather's $\beta$-function, does there exist a diffeomorphism $\psi\colon \T^2 \to \T^2$  isotopic to the identity such that $\psi^* g_2=e^{f} g_1$ for some $f \in C^\infty(\T^2)$? 
\end{enumerate}
\end{remark}

	\smallskip

\section{Main results for Riemannian geodesic flows}

Let $M$ be a closed manifold and let $g_1$ and $g_2$ be two Riemannian metrics on $M$. First of all, we need to introduce some {\it normalization factor} in order to compare two {different metrics and deal with the fact that we can multiply each of them by a constant factor, without changing the dynamical properties}. We define 
\[
\mathcal C_{g_1}(g_2):=\max_{x\in M} \max_{\footnotesize{\substack{v\in T_xM  \\ g_{1,x}(v,v)=1
}}}
 g_{2,x}(v,v)\; = \; \max_{x\in M} \sup_
{\footnotesize{\substack{v\in T_xM  \\ v\neq 0
}}}
\left( \frac{g_{2,x}(v,v)}{g_{1,x}(v,v)}\right),
\]
which can be interpreted as the maximal distortion of the $g_2$-length of tangent vectors at $x$, with respect  to the metric $g_1$. If  $g_1$ and $g_2$ are conformally equivalent metrics on $M$, i.e.,  $g_2 = \varphi \cdot g_1$ for some  positive function $\varphi \in C^{\infty}(M)$, then it is easy to check that
$C_{g_1}(g_2) = \max_{M} \varphi.$\\

We can now state our Main result, which will imply, among other consequences, a pointwise version of Bangert's result. We present our result in terms of Mather's $\beta$-functions, but  it can be equivalently rephrased in terms of the corresponding stable norms by means of~\eqref{eq:beta_stable_relation}.
\\

\begin{mainthm}\label{maintheorem}
Let $M$ be a closed manifold and let $g_1$ and $g_2$ be two Riemannian metrics on $M$. Let $\beta_{g_i}\colon H_1(M;\R)\longrightarrow \R$ be the Mather's $\beta$-function associated to $g_i$,   $i=1,2$. Then:
\begin{enumerate}[label=\textbf{(\roman*)},ref=(\roman*)]
\item\label{un_mthm} We have
\begin{equation}\label{maininequality}
\beta_{g_2}(h) \leq \mathcal C_{g_1}(g_2)\, \beta_{g_1}(h),\quad\forall\,  h \in H_1(M;\R).
\end{equation}
\medskip

\item\label{deux_mthm} Let $h\in H_1(M;\R) \setminus \{0\}$. Then:
{\small \begin{equation*}
\beta_{g_2}(h) = \mathcal C_{g_1}(g_2)\, \beta_{g_1}(h) \quad \Longleftrightarrow\quad
{ \mathfrak M}_{g_1}^h \subseteq  { \mathfrak M}_{g_2}^h 
\;\; \text{and} \;\;
g_{2} = \mathcal C_{g_1}(g_2)\, g_{1} \;\; \text{on}  \;\; \widetilde { \mathcal M}_{g_1}^h \subseteq \widetilde { \mathcal M}_{g_2}^h.
\end{equation*}
}
\medskip
\item\label{trois_mthm} Let us assume that $g_1$ and $g_2$  are conformally equivalent, {\it i.e.}, $g_2=\varphi\cdot \, g_1$ for some positive function $\varphi \in C^\infty (M)$ with
$m := \max_M \varphi$. Then:
$$
\beta_{g_2}(h) \leq m\, \beta_{g_1}(h),\quad \forall\; h\in H_1(M;\R). 
$$
Moreover, for $h\in H_1(M;\R) \setminus \{0\}$ we have
{\small \begin{equation*}
\beta_{g_2}(h) = m\, \beta_{g_1}(h) \quad \Longleftrightarrow \quad
{ \mathfrak M}_{g_1}^h \subseteq  { \mathfrak M}_{g_2}^h 
\;\;  \text{and} \;\;
\varphi \equiv m
\;\;  \text{on}\;\;   {\mathcal M}_{g_1}^h  \subseteq  {\mathcal M}_{g_2}^h.
\end{equation*}
}

\medskip

\end{enumerate}
\end{mainthm}

\medskip

\begin{remark}\label{remarkmainthm}
The Main Theorem implies that $g_1$ and $g_2$ are homothetic (more specifically,  $ g_2 = \mathcal{C}_{g_1}(g_2)\, g_1$) on  the subset of $TM$ (hence, on its closure):
$${\bigcup_{\footnotesize{\substack{h\in H_1(M;\R) \setminus \{0\} \\ \beta_{g_2}(h) = {\mathcal C_{g_1}(g_2)}\, \beta_{g_1}(h) }}}
\!\!\!\!\!\!\!\!\!\!\!\! \R_+ \cdot \widetilde {\mathcal M}_{g_1}^h}, $$
where  we denote $\R_+ \cdot  \widetilde {\mathcal M}_{g_1}^h := \{(x, t v):\;  (x,v)\in  \widetilde {\mathcal M}_{g_1}^h,\; t\in \R_+\}.$

\end{remark}

\bigskip

We can deduce from point~\ref{trois_mthm} of the Main Theorem the following corollary.

\begin{corollary}\label{main-corollary} 
Let $M$ be a closed manifold and let $g_1$ and $g_2$ be  two conformally equivalent Riemannian metrics on $M$, {\it i.e.}, $g_2=\varphi \cdot g_1$ for some positive function $\varphi \in C^\infty (M)$ with $m := \max_M \varphi$.
Let us denote  
$${\mathcal H}_{g_1} := \{ h\in H_1(M;\R)\setminus \{0\}:  \;  {\mathcal M}_{g_1}^h  = M\}.$$
Then,
$\beta_{g_2}(h) = m\, \beta_{g_1}(h)$ for some  $h\in \mathcal H_{g_1}$  
if and only if  $g_1$ and $g_2$ are homothetic, {\it i.e.,} $g_2 = m \, g_1$.\\
\end{corollary}
\medskip

\begin{remark}\label{rem:integrability}
In some cases one can describe the set of rotation vectors ${\mathcal H}_g$ appearing in  Corollary~\ref{main-corollary} or show that it is not empty.  
\begin{enumerate}[label=\textbf{(\roman*)},ref=(\roman*)]
\item\label{un_remque}  In the case of a flat metric $g_0$ on $\T^n$, we have  $\mathcal H_{g_0} := H_1(\T^n;\R)\setminus \{0\}$.
\item\label{deux_remque}  In the case of Liouville metrics on $\T^n$, {\it i.e.}, metrics of the form $g= (f_1(x_1) + \ldots + f_n(x_n))\, (dx_1^2 + \ldots + dx_n^2)$, we have (we identify $H_1(\T^n;\R)) \simeq \R^n$):
$$
\mathcal H_g := \{h\in \R^n: \; h_1\cdot\ldots \cdot h_n \neq 0\}.
$$
\item\label{trois_remque}  The above cases are examples of {\it integrable} metrics (in the sense of Arnold-Liouville). It follows from KAM theory that for sufficiently small perturbations $\tilde g$ of such metrics $g$ the set ${\mathcal H}_{\tilde g}$ is not empty (actually, this set has almost full Lebesgue measure).
\item\label{quatre_remque}  In~\cite{ButlerSorrentino}, the authors considered  Riemannian manifolds $(M,g)$ with a weaker notion of integrability -- which does not require the  manifold to be diffeomorphic to a torus --, under which the set $\mathcal{H}_g$ is  non-empty (actually, it is everything). Consider the following example (see~\cite[Theorem 1.3]{ButlerSorrentino}). Let $G$ be a simply-connected \emph{amenable}  Lie group; recall that a topological group is {amenable} if it admits a left-invariant, finitely additive, Borel probability measure. Let $\Gamma \leq G$ be a lattice subgroup, and $g$ be a metric on $\Gamma \backslash G$ induced by a left-invariant metric on $G$. Then, for all $h \in H_1(\Gamma \backslash G; \mathbb{R})$, the Mather set $\widetilde{\mathcal{M}}^h_{g}$ projects over the whole $\Gamma \backslash G$. Therefore, in this case $\mathcal{H}_g=H_1(\Gamma \backslash G;\R)\setminus \{0\}$.
\end{enumerate}
\end{remark}

\medskip

\begin{remark}
Let us observe that the set ${\mathcal H}_{g}$ is in some sense invariant if we multiply the metric by a constant factor. In fact, let $c>0$ and let $h\in {\mathcal H}_{g}$; then, it is easy to check that $ch  \in {\mathcal H}_{c g}$. Therefore, one can apply Corollary \ref{main-corollary}, by replacing $g_1$ so that $m=\max_M \varphi = 1$.
\end{remark}

\bigskip

On $\T^2$ every Riemannian metric $g$ is conformally equivalent to a flat metric. Let us choose a flat metric $g_0$ in the conformal class of $g$ such that $\mathcal C_{g_0}(g)=1$; this is always possible because a flat metric remains flat under constant rescaling. We call this flat metric $g_0$ the {\it normalized conformally equivalent flat metric} associated to $g$.\\

Main Theorem~\ref{trois_mthm}, Corollary~\ref{main-corollary} and Remark~\ref{rem:integrability}~\ref{un_remque} imply the following result that, in some sense, can be considered a pointwise version  of Bangert's result.

\begin{corollary} \label{maincordim2}
Let $g$ be a Riemannian metric on $\T^2$ and let $g_0$ its normalized conformally equivalent flat metric. Then, $\beta_g(h) \leq \beta_{g_0}(h)$ for every $h\in H_1(\T^2;\R)$. Moreover, equality holds at some  $h\in H_1(\T^2;\R)\setminus \{0\}$ if and only if $g = g_0$, i.e., $g$ is flat.
\end{corollary}

\bigskip

We now prove the Main Theorem.\\

\begin{proof}[\textit{Proof of Main Theorem.}]
\noindent {\bf (i)} If $h=0$, the inequality holds trivially since $\beta_{g_1}(0)=\beta_{g_2}(0)=0$. Let  $h$ be an element in $H_1(M;\R)\setminus \{0\}$ and let $\mu_h$ be a Mather measure for the metric $g_1$ with rotation vector $h$, i.e., such that $\beta_{g_1} (h) = \frac{1}{2} \int_{TM} g_{1,x}(v,v)\, d\mu_h(x,v)$.  Observe that ${\rm supp}\, \mu_h \subset TM$ does not intersect the zero section, since $\widetilde{\mathcal M}_{g_1}^h$ does not intersect the zero section for $h\neq 0$. Then:
\begin{align} \label{ineq_proof1}
\beta_{g_2}(h) &\leq\frac{1}{2} \int_{TM} g_{2,x}(v,v)\, d\mu_h(x,v) \nonumber  \\
&= \frac{1}{2} \int_{TM} \frac{g_{2,x}(v,v)}{g_{1,x}(v,v)} \cdot {g_{1,x}}(v,v)\,  d\mu_h(x,v) \nonumber\\
&\leq \frac{{\mathcal C_{g_1}(g_2)}}{2} \int_{TM}   g_{1,x}(v,v)\, d\mu_h(x,v) \, = \, {\mathcal C_{g_1}(g_2)}\, \beta_{g_1}(h), 
\end{align}
where the first inequality follows from the definition of $\beta_{g_2}(h)$ as the minimum over all closed Borel probability measures of rotation vector $h$.

\smallskip

\noindent  {\bf (ii)} Let $h\in H_1(M;\R) \setminus \{0\}$.\\
$(\Longrightarrow)$ Assume that $\beta_{g_2}(h) = {\mathcal C_{g_1}(g_2)}\, \beta_{g_1}(h)$. Since  $h\neq 0$, we have that $\beta_{g_i}(h) \neq 0$ for $i=1,2$, hence all inequalities in~\eqref{ineq_proof1} are indeed equalities. Then, it follows that on the support of any Mather measure of rotation number $h$ for $g_1$,  the ratio $\frac{g_{2,x}(v,v)}{g_{1,x}(v,v)}$ must be constantly equal to ${\mathcal C_{g_1}(g_2)}$. From the definition of the Mather set of homology class $h$ (recall~\eqref{def:mather_set}) we conclude that 
$$
\frac{g_{2,x}(v,v)}{g_{1,x}(v,v)} \equiv {\mathcal C_{g_1}(g_2)} \qquad \forall\, (x,v) \in \widetilde { \mathcal M}_{g_1}^h.
$$
Observe that the fact that all inequalities in  \eqref{ineq_proof1} are  equalities, also implies that any Mather measure of rotation number $h$ for $g_1$ is also minimizing for $g_2$, {\it i.e.},
$$
{ \mathfrak M}_{g_1}^h \subseteq  { \mathfrak M}_{g_2}^h,
$$
which by~\eqref{def:mather_set} yields  
$$\widetilde { \mathcal M}_{g_1}^h \subseteq \widetilde { \mathcal M}_{g_2}^h.$$
$(\Longleftarrow)$ 
Let  $\mu_h \in { \mathfrak M}_{g_1}^h$. Since ${ \mathfrak M}_{g_1}^h \subseteq  { \mathfrak M}_{g_2}^h$, we have that $\mu_h\in { \mathfrak M}_{g_2}^h$ and therefore: 
\begin{align}
\beta_{g_2}(h) &= \frac{1}{2} \int_{TM} g_{2,x}(v,v)\, d\mu_h(x,v) \nonumber  \\
&= \frac{1}{2} \int_{TM} \frac{g_{2,x}(v,v)}{g_{1,x}(v,v)} \cdot g_{1,x}(v,v)\, d\mu_h(x,v) \nonumber\\
&= \frac{{\mathcal C_{g_1}(g_2)}}{2} \int_{TM}   g_{1,x}(v,v)\, d\mu_h(x,v)\, = \, {\mathcal C_{g_1}(g_2)}\, \beta_{g_1}(h),
\end{align}
where we have used  that $\frac{g_{2,x}(v,v)}{g_{1,x}(v,v)}\equiv {\mathcal C_{g_1}(g_2)} $ on $\widetilde { \mathcal M}_{g_1}^h \subseteq  \widetilde{ \mathcal M}_{g_2}^h$.\\

\noindent  {\bf (iii)} Observe that
if $g_1$ and $g_2$ are conformally equivalent  on $M$, i.e.,  $g_2 = \varphi \cdot g_1$ for some  positive function $\varphi \in C^{\infty}(M)$, then for every $x\in M$ and $v \in T_xM \setminus \{0\}$, we have 
$$
\frac{g_{2,x}(v,v)}{g_{1,x}(v,v)} = \varphi(x)$$
(it is independent of $v$) and in particular $\mathcal{C}_{g_1}(g_2) = \max_{M} \varphi= m$. Since the function $\varphi$ is  independent of $v$, we deduce that  
$$
g_2=\mathcal{C}_{g_1}(g_2)g_1\text{ on }\widetilde { \mathcal M}_{g_1}^h\quad \Longleftrightarrow\quad \varphi\equiv m\text{ on }\mathcal{M}_{g_1}^h.
$$ 
Thus,~\ref{trois_mthm} is a direct consequence of~\ref{un_mthm}-\ref{deux_mthm}.
\qedhere
\end{proof}

\bigskip

\section{Mañé's perturbations of Tonelli Lagrangians} \label{manelagrangian}
In the above discussion, we focused on Riemannian geodesic flows, but the same ideas can be implemented for general Tonelli Lagrangians.
Let $M$ be a closed manifold. A Lagrangian $L\colon TM\, \longrightarrow \,\R$ is said to be {\em Tonelli} if 
\begin{enumerate}
	\item[(i)] $L\in C^2(TM)$;
	\item[(ii)] $L$ is strictly convex in the fibers, i.e., if the second partial  derivative ${\partial^2 L}/{\partial v^2}(x,v)$ is positive definite for all $(x,v)\in TM$;
	\item[(iii)] $L$ is superlinear in each fiber, i.e., $\lim_{\Vert v\Vert \to +\infty}\frac{ L(x,v)}{\Vert v\Vert}=+\infty$, uniformly in $x\in M$.
	\end{enumerate}

The same construction of Mather's $\beta$-function and Mather sets described above for  $L_g$ associated to a Riemannian metric $g$ -- which is clearly Tonelli -- holds for a general Tonelli Lagrangian $L$. 
Let us introduce some notation:
\begin{itemize}
\item we denote by $\beta_L \colon H_1(M;\R) \longrightarrow \R$ the associated $\beta$-function:
\begin{equation}\label{beta gen}
	\beta_g(h):= \min_{\mu\in \rho^{-1}(h)} \int_{TM} L(x,v)\, d\mu\, .
	\end{equation}
\item we denote by  $\mathfrak{M}_{L}^h$ the set of Mather measures of rotation vector $h$, i.e., the set of measures $\mu \in \mathcal{M}(TM)$ realizing the minimum in \eqref{beta gen}.
\item we denote by $ \widetilde{\mathcal{M}}_{L}^h$ the associated Mather set of homology $h$, i.e., the closure of the union of the supports of all Mather measures of rotation vector $h$, and
$\mathcal{M}_L^h := \pi(\widetilde{\mathcal{M}}_{L}^h)$ its projection on $M$.
\end{itemize}
We refer to \cite{SorrLectNotes} for a more detailed description.\\

Given a Tonelli Lagrangian $L$, we can perturb it in the following way. Let $V\in C^2(M)$ and consider the Lagrangian $L_V:= L+V$, which is still of Tonelli type. This kind of perturbation is called {\it in the sense of Ma\~n\'e} (see \cite{Mane}). 
Since we can modify $V$ by a constant without changing the dynamics, we normalize $V$ by assuming that $\max_M V = 0$.

\begin{remark}
Observe that by Maupertuis’ principle, a small perturbation of a geodesic flow on a closed manifold given by a conformal perturbation of the metric $g$ is equivalent to a small perturbation in the Mañé sense of the corresponding Lagrangian $L_g$. In such a framework, given a metric $g$ on $M$, Contreras and Miranda \cite{ConMir} study the structure of the Mather sets for a generic metric in the conformal class of $g$.
\end{remark}

Similarly to what we did above, one can prove the following results.

\begin{theorem}\label{thmmane}
Let $M$ be a closed manifold. Let $L\colon TM\longrightarrow \R$ be a Tonelli Lagrangian and let $V\in C^2(M)$ be such that  $\max_M V = 0$. Then:
\begin{enumerate}[label=\textbf{(\roman*)},ref=(\roman*)]
\item\label{Mane point 1} We have $\beta_{L_V}(h) \leq  \beta_{L}(h)$ for every $h\in H_1(M;\R)$.
\item\label{Mane point 2} Let $h\in H_1(M;\R)\setminus\{0\}$. Then:
{\begin{equation*}
\beta_{L_V}(h) =  \beta_{L}(h) \quad \Longleftrightarrow\quad
{ \mathfrak M}_{L}^h \subseteq  { \mathfrak M}_{L_V}^h 
\;\; \text{and} \;\;
V \equiv 0 \;\; \text{on}  \;\; \widetilde { \mathcal M}_{L}^h \subseteq \widetilde { \mathcal M}_{L_V}^h.
\end{equation*}
}
\item\label{Mane point 3} Let $M=\T^n$ and $L(x,v)=\ell(v)$, i.e., it is a completely integrable system. Let $h\in H_1(M;\R)\setminus\{0\}$. Then, considering $\ell_V:=\ell+V$:
{ \begin{equation*}
\beta_{\ell_V}(h) =  \beta_{\ell}(h) \quad \Longleftrightarrow\quad
V \equiv 0 \;\; \text{on}  \;\; M.
\end{equation*}
}
\end{enumerate}

\end{theorem}

The proof follows closely that of Main Theorem and we omit details. One needs to use the fact that $V\le 0$, which yields $\int_{TM}V\, d\mu\le 0$ and gives the inequality in point~\ref{Mane point 1}. In point~\ref{Mane point 2}, the fact that the previous inequality is actually an equality forces $\int_{TM}V \, d\mu$ to be zero, hence also $V$ (since it is negative). Concerning point~\ref{Mane point 3}, observe that, up to suitably identifying $H_1(\mathbb{T}^n;\mathbb{R})$ with $\mathbb{R}^n$, one get $\beta_{\ell}=\ell$ and $\mathcal M_{\ell}^h=\T^n$ for every $h\in H_1(\T^n;\R)$ (see for instance \cite[Section 4]{SorrUruguay}).\\

\bigskip

\subsection*{Acknowledgements}
This work was largely carried out during a visit by AS to CEREMADE at Université Paris Dauphine-PSL, and he is deeply grateful for the institution's hospitality and stimulating environment.
AS acknowledges the support of the Italian Ministry of University and Research's PRIN 2022 grant ``\textit{Stability in Hamiltonian dynamics and beyond}'' and of the Department of Excellence grant MatMod@TOV (2023-27), awarded to the Department of Mathematics at University of Rome Tor Vergata. AS is member of the INdAM research group GNAMPA and the UMI group DinAmicI. AF is partially supported  by the ANR project CoSyDy (ANR-21-CE40-0014),
the ANR project GALS (ANR-23-CE40-0001) and PEPS project ``{\it Jeuns chercheurs et jeunes chercheuses}'' 2025. AF is member of the UMI group DinAmicI. ML is partially supported  by the ANR project CoSyDy (ANR-21-CE40-0014), the ANR JCJC project PADAWAN (ANR-21-CE40-0012), and by the LESET Math-AMSUD project.

\bibliographystyle{alpha}
\bibliography{BiblioFLS.bib}

\end{document}